\documentclass[12pt]{amsart}

\usepackage{ljm-auth}

\usepackage{amsmath}
\usepackage{amsfonts}
\usepackage{amssymb}
\usepackage{longtable}

\newtheorem{theorem}{Theorem}
\newtheorem{prop}{Proposition}
\newtheorem{tab}{Table}

\def\spa{\mathop\text{{\rm span}}\nolimits}
\def\pr{\mathop\text{\rm pr}\nolimits}
\def\id{\mathop\text{\rm id}\nolimits}

\def\Real{\mathbb{R}}
\def\Co{\mathbb{C}}
\def\H{\mathbb{H}}

\def\g{\mathfrak{g}}
\def\h{\mathfrak{h}}
\def\k{\mathfrak{k}}

\def\z{\mathfrak{z}}

\def\so{\mathfrak{so}}
\def\osp{\mathfrak{osp}}
\def\hosp{\mathfrak{hosp}}
\def\psl{\mathfrak{psl}}
\def\sp{\mathfrak{sp}}
\def\gl{\mathfrak{gl}}
\def\su{\mathfrak{su}}
\def\u{\mathfrak{u}}

\def\pq{\mathfrak{pq}}
\def\pe{\mathfrak{pe}}

\def\spin{\mathfrak{spin}}

\def\p{\partial}

\def\T{{\mathcal T}}
\def\M{{\mathcal M}}

\def\R{{\mathcal R}}
\def\RR{\bar{\mathcal R}}

\author{Anton\,S.\,Galaev}

\crauthor{A.S.\,Galaev}

\tit{Irreducible holonomy algebras of odd Riemannian
supermanifolds}

\shorttit{Irreducible holonomy algebras of odd...}

\begin{document}

\maketit

 \abstract{ Possible irreducible holonomy algebras
$\g\subset\sp(2m,\Real)$ of odd Riemannian supermanifolds and
irreducible subalgebras $\g\subset\gl(n,\Real)$ with non-trivial
first skew-symmetric prolongations are classified. An approach to
the classification of some classes of the holonomy algebras of
Riemannian supermanifolds is discussed. } \notes{0}{

\subclass{58A50, 53C29}%
\keywords{Riemannian supermanifold, Levi-Civita superconnection,
holonomy algebra, Berger superalgebra}%

\thank{Supported by the grant  201/09/P039 of the Grant Agency of
Czech Republic   and by the grant MSM~0021622409 of the Czech
Ministry of Education} }

\section{Introduction}

Berger's classification of  holonomy algebras of Riemannian
manifolds is an important result  that has found applications
 both in geometry and theoretical physics
 \cite{Al1,Ber,Besse,Bryant2,Gibbons09,Joyce07}. For
 pseudo-Riemannian manifolds the corresponding problem is solved
 only in the irreducible case by Berger and in several partial
 cases, e.g. in the Lorentzian signature \cite{ESI}.

Since theoretical physicists discovered supersymmetry,
supermanifolds began to play an important role both in geometry
and physics \cite{DelMor,Leites,Leites1,Manin,Var}. In
\cite{SupHol} the holonomy algebras of linear connections on
supermanifolds are defined. In particular, if $(\M,g)$ is a
Riemannian supermanifold, then its holonomy algebra $\g$ may be
identified with a subalgebra of the orthosymplectic Lie
superalgebra $\osp(p,q|2m)$, where $p+q|2m$ is the dimension of
$\M$ and $(p,q)$ is the signature of the metric $g$ restricted to
the underlying smooth manifold of $\M$. It is natural to pose the
problem of classification of the holonomy algebras
$\g\subset\osp(p,q|2m)$ of Riemannian supermanifolds. Since for
$m=0$ this is the unsolved problem of the differential geometry,
one should consider some restrictions on $\g$. The first natural
restriction is the irreducibility of $\g\subset\osp(p,q|2m)$. Let
us suppose also that
\begin{equation}\label{vidg}\g=(\oplus_i\g_i)\oplus\z,\end{equation}
where $\g_i$ are simple Lie superalgebras of classical type and
$\z$ is a trivial or one-dimensional center (if $m=0$ or $p=q=0$,
then this assumption follows automatically from the irreducibility
of~$\g$).

In the present paper we obtain a classification of possible
irreducible holonomy algebras of odd Riemannian supermanifolds, in
this case $$\g\subset\osp(0|2m)\simeq\sp(2m,\Real)$$ is a usual
Lie algebra. This result is the mirror analog of the Berger
classification. Moreover, the aim of this paper is to collect some
facts that will be needed for the classification of irreducible
holonomy algebras $\g\subset\osp(p,q|2m)$ of the
form~\eqref{vidg}. More precisely, any holonomy algebra
$\g\subset\sp(2m,\Real)$ of an odd Riemannian supermanifold is a
skew-Berger algebra, i.e. $\g$ is spanned by the images of the
space $\bar R(\g)$ that consists of symmetric bilinear  forms on
$\Real^{2m}$ with values in $\g$  satisfying the first Bianchi
identity. These algebras are the analogs of the Berger algebras
that are defined in a similar way \cite{Bryant2,MS99,Sch}. All
previously known irreducible Berger algebras were realized as the
holonomy algebras \cite{Bryant2,Sch}, hence skew-Berger algebras
may be considered as the candidates to the holonomy algebras of
odd supermanifolds. Complex irreducible skew-Berger subalgebras of
$\gl(n,\Co)$ are classified recently in~\cite{Skew-Berger}.

Suppose now that $\g\subset\osp(p,q|2m)$ is of the form
\eqref{vidg} and irreducible. Its even part $$\g_{\bar
0}\subset\so(p,q)\oplus\sp(2m,\Real)$$ preserves the decomposition
$\Real^{p,q}\oplus\Real^{2m}$ into the even and odd parts. For the
most of the representations, $\g_{\bar 0}$ acts diagonally in
$\Real^{p,q}\oplus\Real^{2m}$, i.e. its representations in the
both subspaces are faithful. In \cite{Skew-Berger} it is explained
that in this case $\pr_{\so(p,q)}\g_{\bar 0}\subset\so(p,q)$ is a
Berger algebra and $\pr_{\sp(2m,\Real)}\g_{\bar
0}\subset\sp(2m,\Real)$ is a skew-Berger algebra. Since $\g_{\bar
0}$ is reductive, $\pr_{\so(p,q)}\g_{\bar 0}$ is known.  Below we
will show that to know all possible $\pr_{\sp(2m,\Real)}\g_{\bar
0}$, it is enough to classify irreducible skew-Berger subalgebras
 $\h\subset\sp(2k,\Real)$ and to classify skew-Berger subalgebras
$\h\subset\sp(2m,\Real)$ that preserve a decomposition
$\Real^{2m}=W\oplus W_1$ into the direct sum of two Lagrangian
subspaces and act diagonally in $W\oplus W_1$. The last
classification can be reduced to the classification of irreducible
subalgebras $\h\subset\gl(n,\Real)$ with the non-trivial first
skew-symmetric  prolongation
$$\h^{[1]}=\{\varphi\in(\mathbb{R}^n)^*\otimes\h|\varphi(x)y=-\varphi(y)x\text{
for all } x,y\in\mathbb{R}^n\}.$$ For subalgebras
$\h\subset\so(n,\Real)$ the first skew-symmetric  prolongation is
studied in \cite{Nagy}, where some applications are obtained.
Thus, knowing $\pr_{\so(p,q)}\g_{\bar 0}$ and
$\pr_{\sp(2m,\Real)}\g_{\bar 0}$, and using the theory of
representations of simple Lie superalgebras \cite{KacRepr}, it is
possible to find $\g\subset\osp(p,q|2m)$. These ideas and the
results of this paper will allow to classify irreducible
subalgebras $\g\subset\osp(p,q|2m)$ of the form \eqref{vidg}.

The paper has the following structure. In Section \ref{SecPrel} we
give necessary preliminaries.  Section~\ref{secSymspaces} deals
with  odd Riemannian symmetric superspaces. In Section~\ref{Prol}
the classification of irreducible subalgebras
$\g\subset\gl(n,\Real)$ with non-trivial first skew-symmetric
prolongation is obtained. In Section~\ref{SBa}  irreducible not
symmetric skew-Berger subalgebras $\g\subset\sp(2m,\Real)$, i.e.
possible irreducible holonomy algebras of not symmetric odd
Riemannian supermanifolds, are classified.

\section{Preliminaries}\label{SecPrel}

{\bf Odd Riemannian supermanifolds, connections, holonomy
algebras.} First we rewrite some  general definitions from the
theory of supermanifolds \cite{Leites,Leites1,Manin} for the case
of odd supermanifolds. A connected supermanifold $\M$ of dimension
$0|k$ is a pair $$(\{x\},\Lambda(k))$$ where $x$ is the only point
of the manifold and $\Lambda(k)$ is a Grassman superalgebra of $k$
generators, which is considered as the superalgebra of functions
on $\M$. Such supermanifolds $\M$ are called odd supermanifolds.
If $\xi^1,\dots,\xi^k$ are generators of $\Lambda(k)$, then
$$\Lambda(k)=\oplus_{i=0}^k\Lambda^i\Real^k=\Lambda(k)_{\bar 0}\oplus\Lambda(k)_{\bar 1},$$
where $\Real^k=\spa\{\xi^1,\dots,\xi^k\}$, and $\Lambda(k)_{\bar
0}$ and $\Lambda(k)_{\bar 1}$ are spanned by the elements of even
and odd degree, respectively. The elements of $\Lambda(k)_{\bar
0}$ and $\Lambda(k)_{\bar 1}$ are called homogeneous. For a
homogeneous element $f$, the parity $|f|\in\mathbb{Z}_2=\{\bar
0,\bar 1\}$ is defined to be $|\bar 0|$ or $|\bar 1|$ if
$f\in\Lambda(k)_{\bar 0}$ or $f\in\Lambda(k)_{\bar
1}\backslash\{0\}$, respectively. It holds
$|\xi^1|=\cdots=|\xi^k|=\bar 1$ and
$$fh=(-1)^{|f||h|}hf$$ for all homogenous $f,h\in\Lambda(k)$. In particular, $\xi^i\xi^j=-\xi^j\xi^i$.
Any function $f\in\Lambda(k)$ can be written in the form
$$f=\sum_{r=0}^m\sum_{\alpha_1<\cdots<\alpha_r}f_{\alpha_1\dots\alpha_r}\xi^{\alpha_1}\cdots\xi^{\alpha_r},$$ where
 $f_{\alpha_1\dots\alpha_r}\in\Real$. By definition,  the value of $f$ at the point $x$ is $f(x)=f_\emptyset$.
The functions $\xi^1,\dots,\xi^k$ are called coordinates on $\M$.
The $\Lambda(k)$-supermodule $$\T_\M=(\T_\M)_{\bar
0}\oplus(\T_\M)_{\bar 1}$$ of vector fields on $\M$ consists of
$\Real$-linear maps $X:\Lambda(k)\to\Lambda(k)$ such that the
homogeneous $X$ satisfy
$$X(fh)=(Xf)h+(-1)^{|X||f|}f(Xh)$$ for all homogeneous functions $f,g$.
The vector fields $\p_{\xi^1},\dots,\p_{\xi^k}$ are defined in the
obvious way. These vector fields are odd. The tangent space
$T_x\M$ to $\M$ at the point $x$ can be identified with
$\spa_\Real\{\p_{\xi^1},\dots,\p_{\xi^k}\}$ and it is an odd
vector superspace. It holds $$\T_\M=\Lambda(k)\otimes T_x\M.$$ The
value of a vector field $X=f_\alpha\p_{\xi^\alpha}$ at the point
$x$ is defined as $X_x=f_\alpha(x)(\p_{\xi^\alpha})_x\in T_x\M$.

{\it A connection} on $\M$ is an even $\Real$-linear map
$$\nabla:\T_\M\otimes_\Real \T_\M\to\T_\M$$ of
$\Real$-supermodules such that
$$\nabla_{fY}X=f\nabla_YX\quad\text{and}\quad\nabla_YfX=(Yf)X+(-1)^{|Y||f|}f\nabla_YX$$
for all homogeneous functions $f$ and vector fields $X$, $Y$ on
$\M$. The curvature tensor $R$ of $\nabla$ and its covariant
derivatives $\nabla^rR$ and their values at the point $x$ are
defined in the usual way. From \cite{SupHol} it follows that {\it
the holonomy algebra} $\g$ of the connection $\nabla$ at the point
$x$ can be defined in the following way:
\begin{multline*}\g=\spa\left\{\nabla^r_{\p_{\xi^{\alpha_r}},\dots,\p_{\xi^{\alpha_1}}}R_x(\p_{\xi^\beta},\p_{\xi^\gamma})\left|
\begin{matrix}\,0\leq r\leq k,\,1\leq\beta,\gamma\leq
k,\\
 1\leq \alpha_1<\dots<\alpha_r\leq
 k\end{matrix}\right\}\right.\\ \subset\gl(T_x\M)\simeq\gl(0|k,\Real).\end{multline*}
Thus $\g$ is a usual Lie algebra acting in the odd superspace
$T_x\M$. Considering the isomorphism $\Pi T_x\M\simeq \Real^k$, we
get $\g\subset\gl(k,\Real)$. Here $\Pi$ is the parity changing
functor. The holonomy group of the connection $\nabla$ at the
point $x$ is defined as the corresponding connected Lie subgroup
of ${\rm Gl}(k,\Real)$.

{\it A Riemannian supermetric} on $\M$ is an even linear map
$$g:\odot^2\T_\M\to \Lambda(k)$$ such that its value
$$\omega=g_x\in\odot^2T^*_x\M$$ is non-degenerate. Since $T_x\M$
is an odd vector superspace, $\omega$ is a symplectic form on $\Pi
T_x\M\simeq\Real^k$. Hence in this case $k$ must be even, $k=2m$.
On such Riemannian supermanifold there exists the Levi-Civita
superconnection $\nabla$. For its holonomy algebra it holds
$$\g\subset\osp(0|2m,\Real)\simeq\sp(2m,\Real).$$

{\bf Skew-Berger algebras.} The main task of this paper is to
classify possible irreducible holonomy algebras
$\g\subset\sp(2m,\Real)$ of odd Riemannian supermanifolds. This
can be done using the following algebraic properties of the
representation $\g\subset\sp(2m,\Real)$.

 Let $V$ be a real or complex vector space and $\g\subset\gl(V)$ a subalgebra.
The space of skew-symmetric curvature tensors of type $\g$ is
defined as follows: $$\bar\R(\g)=\left\{R\in\odot^2 V^*\otimes
\g\,\left|\begin{array}{c}R(X, Y)Z+R(Y, Z)X+R(Z, X)Y=0\\ \text{
}\text{ for all  } X,Y,Z\in V\end{array}\right\}\right..$$ The
subalgebra $\g\subset\gl(V)$ is called {\it a skew-Berger
subalgebra} if it is spanned by the images of the elements
$R\in\bar\R(\g)$. The space $\R(\g)$ and the notion of the Berger
algebra (or more generally of a Berger superalgebra) is defined in
the same way, the only difference is that $R$  is a (super)
skew-symmetric bilinear form with values in $\g$.
 Obviously $\bar\R(\g)=\R(\g\subset\gl(\Pi V))$ and $\g\subset\gl(V)$ is a
skew-Berger algebra if and only if $\g\subset\gl(\Pi V)$ is a
Berger superalgebra.

Let $\g\subset\gl(k,\Real)$ Consider the space of linear maps from
$\Real^k$ to $\RR(\g)$ satisfying the second Bianchi identity,
$$\RR^\nabla(\g)=\left\{S\in  (\Real^k)^* \otimes\RR(\g)\left|\begin{matrix}S_X(Y,Z)+S_Y(Z,X)+S_Z(X,Y)=0\\\text{for all  } X,Y,Z\in \Real^k
\end{matrix}\right\}\right..$$ A skew-Berger subalgebra
$\g\subset\gl(n,\Real)$ is called {\it symmetric}  if
$\RR^\nabla(\g)=0$. Note that if an odd supermanifold $\M$ is
endowed with a torsion-free connection $\nabla$ and $\g$ is its
holonomy algebra, then
$\nabla^r_{\p_{\xi^{\alpha_r}},\dots,\p_{\xi^{\alpha_1}}}R_x\in\RR(\g)$
and
$\nabla^r_{\p_{\xi^{\alpha_r}},\dots,\p_{\xi^{\alpha_2}},\cdot}R_x\in
\RR^\nabla(\g)$. In particular, $\g\subset\gl(k,\Real)$ is a
skew-Berger algebra.

Let $\omega$ be the standard symplectic form on $\Real^{2m}$. A
subalgebra $\g\subset\sp(2m,\Real)$ is called {\it
weakly-irreducible} if it does not preserve any proper
non-degenerate subspace of $\Real^{2m}$. The next theorem is the
partial case of the Wu Theorem for Riemannian supermanifolds
proved in \cite{SupHol}.

\begin{theorem}  Let $\g\subset\sp(2m,\Real)=\sp(V)$ be a skew-Berger subalgebra, then there is a decomposition $$V=V_0\oplus
V_1\oplus\cdots\oplus V_r$$ into a direct sum of symplectic
subspaces and a decomposition $$\g=\g_1\oplus\cdots\oplus \g_r$$
into a direct sum of ideals such that $\g_i$ annihilates $V_j$ if
$i\neq j$ and $\g_i\subset \sp(V_i)$ is a weakly-irreducible
Berger subalgebra.

If $\g\subset\gl(k,\Real)$ is the holonomy algebra of an odd
Riemannian supermanifold $(\M,g)$, then the above decompositions
define a decomposition of $(\M,g)$ into the product of a flat odd
Riemannian supermanifold and of odd Riemannian supermanifolds with
the weakly-irreducible holonomy algebras~$\g_i\subset\sp(V_i)$.
\end{theorem}

Let $\g\subset\sp(2m,\Real)=\sp(V)$ be a skew-Berger subalgebra.
By the above theorem, we may assume that it is
 weakly-irreducible. Suppose that it is not
irreducible. Suppose also that $\g$ is a reductive Lie algebra.
Since $\g$ is not irreducible,  it preserves a degenerate subspace
$W\subset V$. Consequently, $\g$ preserves the isotropic subspace
$L=W\cap W^\perp$ ($W^\perp$ is defined using $\omega$). Since
$\g$ is totally reducible, there exists a complementary invariant
subspace $L'\subset V$. Since $\g$ is weakly-irreducible, the
subspace $L'$ is degenerate. If $L'$ is not isotropic, then $\g$
preserves the kernel of the restriction of $\omega$ to $L'$ and
$\g$ preserves a complementary subspace in $L'$ to this kernel,
which is non-degenerate. Hence $L'$ is isotropic and $V=L\oplus
L'$ is the direct of two Lagrangian subspaces. The form $\omega$
on $V$ allows to identify $L'$ with the dual space $L^*$ and the
representations of $\g$ on $L$ and $L'$ are dual. Since
$\g\subset\sp(V)$ is weakly-irreducible, the representation
$\g\subset\gl(L)$ is irreducible. Let $R\in\RR(\g)$. From the
Bianchi identity it follows that $R(x,y)=0$ and
$R(\varphi,\psi)=0$ for all $x,y\in L$ and $\varphi,\psi\in L^*$.
Moreover, for each fixed $\varphi\in L^*$ it holds
$R(\cdot,\varphi)\in(\g\subset\gl(L))^{[1]}$, where
$(\g\subset\gl(L))^{[1]}$ is the first skew-symmetric prolongation
for the representation $\g\subset\gl(L)$ (similarly, for each
fixed $x\in L$ it holds $R(\cdot,x)\in(\g\subset\gl(L^*))^{[1]}$).
Consequently, $(\g\subset\gl(L))^{[1]}\neq 0$ and such algebras
are classified in Section \ref{Prol}. Thus we will get
classification of all reductive skew-Berger subalgebras
$\g\subset\sp(2m,\Real)$.

\section{Odd symmetric superspaces and simple Lie
superalgebras}\label{secSymspaces}

Symmetric superspaces are studied in \cite{Serganova,Leites1,Go}.
A class of odd symmetric superspaces is considered in
\cite{Cortes2}.

An odd supermanifold $(\M,\nabla)$ is called symmetric if $\nabla
R=0$. In this case for the holonomy algebra $\g$ we have
$$\g=\spa\{R_x(\p_{\xi^\beta},\p_{\xi^\gamma})| 1\leq \beta,\gamma\leq k\}=R_x(T_x\M,T_x\M).$$
Moreover, $\g$ annihilates $R_x\in\R(\g)$ \cite{SupHol}.

Similarly  as in \cite{Sch} it can be shown that {\it if
$\g\subset\gl(k,\Real)$ is irreducible skew-Berger algebra such
that the representation of $\g$ in $\R(\g)$ is trivial, then $\g$
is a symmetric skew-Berger algebra.}

From \cite{SupHol} it follows that {\it if the holonomy algebra
$\g$ of a torsion-free connection $\nabla$ on an odd supermanifold
$\M$ is a symmetric skew-Berger algebra, then $(\M,\nabla)$ is
symmetric.}

Let $(M,\nabla)$ be an odd symmetric  supermanifold. Define the
Lie superalgebra
$$\k=\g\oplus\Pi\Real^k$$ with the Lie superbrackets $$[\Pi X,\Pi
Y]=R(X,Y),\quad [A,\Pi X]=\Pi(AX),\quad [A,B]=[A,B]_\g,$$ where
$A,B\in\g$ and $X,Y\in\Real^k$. Conversely, let $\k$ be a Lie
superalgebra with the even part $\g$ and the odd part
$\Pi\Real^k$.  Let $G$ be the connected Lie group with the Lie
algebra $\g$ and  ${\mathcal K}$ be the connected Lie supergroup
with the Lie superalgebra $\k$. The Lie supergroup $\mathcal K$
can be given by the Harish-Chandra pair $(G,\k)$ \cite{Go}. The
factor superspace $\M={\mathcal K}/G$ is an odd supermanifold and
it admits a  unique symmetric superconnection~\cite{Leites1}.

Thus we obtain a  one to one correspondence between connected odd
symmetric superspaces $(\M,\nabla)$ and Lie superalgebras
$\k=\g\oplus\Pi\Real^k$. Moreover,  $\Pi\Real^k$ is the tangent
space to $\M$ and $\g\subset\gl(k,\Real)$ is the holonomy algebra.
The space $(\M,\nabla)$ is Riemannian if and only if
$\g\subset\sp(k,\Real)$.

Let $\k=\k_{\bar 0}\oplus \k_{\bar 1}$ be a (real or complex)
simple Lie superalgebra. It is of {\it classical type} if the
representation of $\k_{\bar 0}$ on $\k_{\bar 1}$ is totally
reducible. In this case, $\k$ is of {\it type I}  if the
representation of $\k_{\bar 0}$ on $\k_{\bar 1}$ is irreducible;
$\k$ is of {\it type II} if $\k_{\bar 0}$ preserves a
decomposition $$\k_{\bar 1}=\k_{-1}\oplus\k_1$$ such that the
representations of $\k_{\bar 0}$ in  $\k_{-1}$ and $\k_{1}$ are
faithful and irreducible.

Let $\g$ be a reductive Lie algebra. Suppose that
$\g\subset\gl(k,\Real)$ is irreducible, or there exists a
$\g$-invariant decomposition  $$\Pi\Real^k=\k_{-1}\oplus\k_{1}$$
such that the representations of $\g$ on  $\k_{-1}$ and $\k_{1}$
are faithful and irreducible. Suppose that there exists
$R\in\bar\R(g)$ such that $\g$ annihilates $R$ and
$R(\Real^k,\Real^k)=\g$, then the Lie superalgebra
$\k=\g\oplus\Pi\Real^k$ defined as above  is simple, this follows
from Propositions 1.2.7 and 1.2.8 from \cite{Kac}.

Thus we have reduced the classification of weakly-irreducible
reductive subalgebras $\g\subset\gl(k,\Real)$ admitting elements
$R\in\bar\R(g)$ such that $\g$ annihilates $R$ and
$R(\Real^k,\Real^k)=\g$ to the classification of real simple Lie
superalgebras of classical type. Remark that we are interested
only in the case $\g\subset\sp(k,\Real)$. If
$\g\subset\gl(k,\Real)$ is not irreducible, then
$\g\subset\sp(k,\Real)$ if and only if $\k_{-1}\simeq \k_1^*$.

Real simple Lie superalgebras of classical type are exhausted by
simple complex Lie superalgebras of classical type considered as
real Lie superalgebras and by real forms of complex Lie
superalgebras of classical type. The real forms are found in
\cite{Parker80}. To make the exposition complete we list these
algebras in Tables 1--4.

{\footnotesize

\begin{tab} Simple complex Lie superalgebras of type I\\
\begin{tabular}{|c|c|c|c|} \hline $\g$&$\g_{\bar 0}$&$\g_{\bar
1}$&restriction\\ \hline
$\osp(n|2m,\Co)$&$\so(n,\Co)\oplus\sp(2m,\Co)$&$\Co^n\otimes\Co^{2m}$&$n\neq
2$\\
\hline
$\osp(4|2,\alpha,\Co)$&$\mathfrak{sl}(2,\Co)\oplus\mathfrak{sl}(2,\Co)\oplus\mathfrak{sl}(2,\Co)$&$\Co^2\otimes\Co^{2}\otimes\Co^{2}$&
$\alpha\in\Co\backslash\{0,-1\}$\\ \hline
$F(4)$&$\so(7,\Co)\oplus\mathfrak{sl}(2,\Co)$&$\Co^8\otimes\Co^2$&\\
\hline
$G(3)$&$G_2\oplus\mathfrak{sl}(2,\Co)$&$\Co^7\otimes\Co^2$&\\
\hline
$\pq(n,\Co)$&$\mathfrak{sl}(n,\Co)$&$\mathfrak{sl}(n,\Co)$&$n\geq
3$
\\\hline\end{tabular}\end{tab}

\begin{tab} Simple complex Lie superalgebras of type II\\
\begin{tabular}{|c|c|c|c|c|}\hline $\g$&$\g_{\bar
0}$&$\g_1$&$\g_{-1}$&restriction\\ \hline
$\mathfrak{sl}(n|m,\Co)$&$\mathfrak{sl}(n,\Co)\oplus
\mathfrak{sl}(m,\Co)\oplus \Co$&$\Co^{n}\otimes\Co^{m*}$&$\Co^{n*}
\otimes\Co^{m}$&$n\neq m$\\ \hline
$\psl(n|n,\Co)$&$\mathfrak{sl}(n,\Co)\oplus
\mathfrak{sl}(n,\Co)$&$\Co^{n}\otimes\Co^{n*}$&$\Co^{n*}\otimes\Co^{n}$&\\
\hline
$\osp(2|2m,\Co)$&$\so(2,\Co)\oplus\sp(2m,\Co)$&$\Co^{2m}$&$\Co^{2m}$&\\
\hline$\pe(n,\Co)$&$\mathfrak{sl}(n,\Co)$&$\odot^2\mathfrak{sl}(n,\Co)$&$\Lambda^2\mathfrak{sl}(n,\Co)^*$&$n\geq
3$\\\hline
\end{tabular}\end{tab}

\begin{tab} Simple real Lie superalgebras of type I\\
\begin{tabular}{|c|c|c|c|} \hline $\g\otimes\Co$& $\g$&$\g_{\bar 0}$&$\g_{\bar
1}$\\ \hline $\osp(n|2m,\Co)$& $\osp(r,n-r|2m,\Real)$&
$\so(r,n-r)\oplus
\sp(2m,\Real)$&$\Real^{r,n-r}\otimes\Real^{2m}$\\\hline
$\osp(2n|2m,\Co)$& $\hosp(r,m-r|n)$& $\so(n,\mathbb{H})\oplus
\sp(r,m-r)$&$\mathbb{H}^{n}\otimes_\H\mathbb{H}^{r,m-r}$\\\hline
$\osp(1|2m,\Co)$&$\osp(1|2m,\Real)$&$\sp(2m,\Real)$&$\Real^{2m}$\\\hline
$\osp(4|2,\alpha,\Co)$&&$\mathfrak{sl}(2,\Real)\oplus\mathfrak{sl}(2,\Real)\oplus\mathfrak{sl}(2,\Real)$&$\Real^2\otimes\Real^2\otimes\Real^2$\\
&&$\su(2)\oplus\su(2)\oplus\mathfrak{sl}(2,\Real)$&$\Real^4\otimes\Real^2$\\
&&$\mathfrak{sl}(2,\Co)\oplus\mathfrak{sl}(2,\Real)$&$\Real^{1,3}\otimes\Real^2$\\\hline
$F(4)$&&$\mathfrak{sl}(2,\Real)\oplus\so(7)$&$\Real^2\otimes\Delta_7$\\
&&$\mathfrak{sl}(2,\Real)\oplus\so(3,4)$&$\Real^2\otimes\Delta_{3,4}$\\
&&$\su(2)\oplus\so(2,5)$&$\Co^2\otimes\Delta_{2,5}$\\
&&$\su(2)\oplus\so(1,6)$&$\Co^2\otimes\Delta_{1,6}$\\\hline
$G(3)$&&$\mathfrak{sl}(2,\Real)\oplus G_2$&$\Real^2\otimes\Real^7$\\
&&$\mathfrak{sl}(2,\Real)\oplus
G^*_{2(2)}$&$\Real^2\otimes\Real^{3,4}$\\\hline
$\pq(n,\Co)$&$\pq(n,\Real)$&$\mathfrak{sl}(n,\Real)$&$\mathfrak{sl}(n,\Real)$\\
&&$\su(p,n-p)$&$\su(p,n-p)$\\
&$\pq(n,\H)$&$\mathfrak{sl}(\frac{n}{2},\mathbb{H})$&$\mathfrak{sl}(\frac{n}{2},\mathbb{H})$\\\hline
$\mathfrak{sl}(n|m,\Co)$&$\su(s,n-s|r,m-r)$&$\su(s,n-s)\oplus\su(r,m-r)\oplus\Real
i$&$\Co^{s,n-s}\otimes\Co^{r,m-r}$\\\hline
$\mathfrak{psl}(n|n,\Co)$&$\mathfrak{psu}(s,n-s|r,n-r)$&$\su(s,n-s)\oplus\su(r,n-r)$&$\Co^{s,n-s}\otimes\Co^{r,n-r}$\\
&&$\mathfrak{sl}(n,\Co)$&$\Co^n\otimes\overline{\Co^{n}}$\\\hline
$\osp(2|2m,\Co)$&$\hosp(r,m-r|1)$&$\Real
i\oplus\sp(r,m-r)$&$\mathbb{H}^{r,m-r}$\\\hline

\end{tabular}\end{tab}

\begin{tab} Simple real Lie superalgebras of type II\\
\begin{tabular}{|c|c|c|c|} \hline $\g\otimes\Co$& $\g$&$\g_{\bar 0}$&$\g_{\bar
1}$\\ \hline
$\mathfrak{sl}(n|m,\Co)$& $\mathfrak{sl}(n|m,\Real)$& $\mathfrak{sl}(n,\Real)\oplus \mathfrak{sl}(m,\Real)\oplus \Real$&$\Real^{n}\otimes\Real^{m*}\oplus \Real^{n*}\otimes\Real^{m}$\\
&$\mathfrak{sl}(\frac{n}{2}|\frac{m}{2},\mathbb{H})$&
$\mathfrak{sl}(\frac{n}{2},\H)\oplus
\mathfrak{sl}(\frac{m}{2},\H)\oplus
\Real$&$\mathbb{H}^\frac{n}{2}\otimes_\H\mathbb{H}^{\frac{m}{2}*}\oplus\mathbb{H}^{\frac{n}{2}*}\otimes_\H\mathbb{H}^{\frac{m}{2}}$\\\hline
$\osp(2|2m,\Co)$&$\osp(2|2m,\Real)$&$\Real\oplus\sp(2m,\Real)$&$\Real^{2m}\oplus\Real^{2m*}$\\\hline
$\pe(n,\Co)$&$\pe(n,\Real)$&$\mathfrak{sl}(n,\Real)$&$\odot^2\mathfrak{sl}(n,\Real)\oplus\Lambda^2\mathfrak{sl}(n,\Real)^*$\\
&$\pe(\frac{n}{2},\Real)$&$\mathfrak{sl}(\frac{n}{2},\H)$&$\odot^2\mathfrak{sl}(\frac{n}{2},\H)\oplus\Lambda^2\mathfrak{sl}(\frac{n}{2},\H)^*$\\\hline
\end{tabular}\end{tab}

}

\section{Skew-symmetric prolongations of Lie algebras}\label{Prol}

Irreducible subalgebras $\g\subset\gl(n,\mathbb{F})$
($\mathbb{F}=\Real$ or $\Co$) with non-trivial prolongations
$$\g^{(1)}=\{\varphi\in(\mathbb{F}^n)^*\otimes\g|\varphi(x)y=\varphi(y)x\text{
for all } x,y\in\mathbb{F}^n\}$$ are well known, see e.g.
\cite{Bryant2}. Here we classify irreducible subalgebras
$\g\subset\gl(n,\mathbb{F})$
 such that the skew-symmetric
prolongation
$$\g^{[1]}=\{\varphi\in(\mathbb{F}^n)^*\otimes\g|\varphi(x)y=-\varphi(y)x\text{
for all } x,y\in\mathbb{F}^n\}$$ of $\g$ is non-zero.

Irreducible subalgebras $\g\subset\so(n,\Real)$ with non-zero skew-symmetric prolongations are classified in \cite{Nagy}.
These subalgebras are exhausted by the whole orthogonal Lie  algebra $\so(n,\Real)$ and by the adjoint representations of
compact simple Lie algebras. 

Irreducible subalgebras $\g\subset\gl(n,\mathbb{C})$ with
$\g^{[1]}\neq 0$ are classified in \cite{Skew-Berger}. We give
this list in Table \ref{sProlC}. To get this result it was used
that $\g^{[1]}$ coincides with $\Pi(\g\subset\gl(0|n,\Co))^{(1)}$
and the fact that the whole Cartan prolong $$\g_*=\Pi
V\oplus\g\oplus (\g\subset\gl(0|n,\Co))^{(1)}\oplus
(\g\subset\gl(0|n,\Co))^{(2)}\oplus\cdots$$ is an irreducible
transitive Lie superalgebra with the consistent
$\mathbb{Z}$-grading and $\g_{1}\neq 0$. All such
$\mathbb{Z}$-graded Lie superalgebras are classified in
\cite{Kac}. The second prolongation $\g^{[2]}$ is defined in the
obvious way.

{\footnotesize

\begin{tab}\label{sProlC}
Complex irreducible subalgebras $\g\subset\gl(V)$ with
$\g^{[1]}\neq 0$.\\
$
\begin{array}{|c|cc|c|c|}\hline 
\g& V&&\g^{[1]}&\g^{[2]}
\\\hline
\mathfrak{sl}(n,\Co)&\Co^n,&n\geq
3&(\Co^n\otimes\Lambda^2(\Co^n)^*)_0&(\Co^n\otimes\Lambda^3(\Co^n)^*)_0
\\\hline
\gl(n,\Co)&\Co^n,&n\geq 2
&\Co^n\otimes\Lambda^2(\Co^n)^*&\Co^n\otimes\Lambda^3(\Co^n)^*
\\\hline
\mathfrak{sl}(n,\Co)&\odot^2\Co^n,&n\geq 3 &\Lambda^2(\Co^n)^*&0
\\\hline
\gl(n,\Co)&\odot^2\Co^n,&n\geq 3
&\Lambda^2(\Co^n)^*&0
\\\hline
\mathfrak{sl}(n,\Co)&\Lambda^2\Co^n,&n\geq 5&\odot^2(\Co^n)^*&0
\\\hline
\gl(n,\Co)&\Lambda^2\Co^n,&n\geq 5&\odot^2(\Co^n)^*&0
\\\hline
\mathfrak{sl}(n,\Co)\oplus\mathfrak{sl}(m,\Co)\oplus\Co&\Co^n\otimes\Co^m,&
n,m\geq 2 &V^*&0\\ &&n\neq m&&
\\\hline
\mathfrak{sl}(n,\Co)\oplus\mathfrak{sl}(n,\Co)&\Co^n\otimes\Co^n,&
n\geq 3&V^*&0
\\\hline
\mathfrak{sl}(n,\Co)\oplus\mathfrak{sl}(n,\Co)\oplus\Co&\Co^n\otimes\Co^n,&n\geq
3&V^*&0
\\\hline
\so(n,\Co)&\Co^n,&n\geq 4&\Lambda^3 V^*&\Lambda^4 V^*
\\\hline
\so(n,\Co)\oplus\Co&\Co^n,&n\geq 4&\Lambda^3 V^*&\Lambda^4 V^*
\\\hline
\sp(2n,\Co)\oplus\Co&\Co^{2n}&n\geq 2&V^*&0
\\\hline
\g \text{ is simple }&\g& &\Co\id&0
\\\hline
\g\oplus\Co,\,\,\g \text{ is simple }&\g& &\Co\id&0
\\\hline\end{array} $
\end{tab}
}

Now we classify irreducible subalgebras
$\g\subset\gl(n,\mathbb{R})$ with $\g^{[1]}\neq 0$. Let
$\g\subset\gl(n,\mathbb{R})$ be such subalgebra. If this
representation is absolutely irreducible, i.e. $\Real^n$ does not
admit a complex structure commuting with the elements of $\g$,
then $\g\otimes\Co\subset\gl(n,\mathbb{C})$ is an irreducible
subalgebra, and $(\g\otimes\Co)^{[1]}=\g^{[1]}\otimes\Co\neq 0$.
Note that if in this case the representation
$\g\subset\gl(n,\mathbb{R})$ is different from the adjoint one and
from the standard representation of $\so(n,\Real)$, then
$(\g\otimes\Co)^{(1)}\neq 0$ or $(\g\otimes\Co\oplus\Co)^{(1)}\neq
0$ (then $\g^{(1)}\neq 0$). Hence absolutely irreducible
subalgebras $\g\subset\gl(n,\mathbb{R})$ with $\g^{[1]}\neq 0$ are
exhausted by absolutely irreducible subalgebras
$\g\subset\gl(n,\mathbb{R})$ (up to the center of $\g$) with
$\g^{[1]}\neq 0$, by the adjoint representations of real forms of
complex simple Lie algebras, and by
$\so(p,n-p)\subset\gl(n,\Real)$. The result is given in Table
\ref{sProlR}, where we use the following notation from
\cite{Bryant2}:
$$H_n(\Co)=\{A\in Mat_n(\Co)|A^*=A\},\quad S_n(\mathbb{H})=\{A\in Mat_n(\mathbb{H})|A^*=-A\},$$
$$A_n(\mathbb{H})=\{A\in Mat_n(\mathbb{H})|A^*=A\}.$$ The first and the second skew-symmetric prolongations  can be found from the relation $(\g\otimes\Co)^{[k]}=\g^{[k]}\otimes\Co$.

Suppose that the representation $\g\subset\gl(n,\mathbb{R})$  is
non-absolutely irreducible, i.e. $E=\Real^n$  admits a complex
structure $J$ commuting with the elements of $\g$. In this case
the complexificated space $E\otimes\Co$ admits the decomposition
$E\otimes\Co=V\oplus\bar V$, where $V$ and $\bar V$ are the
eigenspaces of the extension of $J$ to $E\otimes \Co$
corresponding to the eigenvalues $i$ and $-i$, respectively. The
Lie algebra $\g\otimes \Co$ preserves this decomposition. Consider
the ideal $\g_1=\g\cap J\g\subset\g$. Since $\g$ is reductive,
there is an ideal $\g_2\subset\g$ such that $\g=\g_1\oplus\g_2$.
The Lie algebra $\g_1\otimes\Co$ admits the decomposition
$\g_1\otimes\Co=\g'_1\oplus\g''_1$ into the eigenspaces of the
extension of $J$ to $\g_1\otimes\Co$ corresponding to the
eigenvalues $i$ and $-i$, respectively. It is easy to see that
$\g'_1$ annihilates $V$, $\g''_1$ annihilates $\bar V$, and
$\g_2\otimes \Co$ acts diagonally in $V\oplus\bar V$. We
immediately conclude that
$$(\g\otimes\Co)^{[1]}=(\g'_1\subset\gl(\bar V))^{[1]}\oplus
(\g''_1\subset\gl(V))^{[1]}.$$ It is clear that the representation
of $\g'_1\oplus(\g_2\otimes\Co)$ in $\bar V$ is irreducible. If
$\dim\g_2\geq 2$, then this representation is of the form of the
tensor product of irreducible representations of   $\g'_1$ and
$\g_2\otimes\Co$. Obviously, in this case $(\g'_1\subset\gl(\bar
V))^{[1]}=0$, similarly $(\g''_1\subset\gl(V))^{[1]}=0$. We
conclude that $\dim\g_2\leq 1$, and
$\g_1\subset\gl(\frac{n}{2},\mathbb{C})\subset \gl(n,\mathbb{R})$
is a complex subalgebra considered as the real one.

Thus irreducible subalgebras $\g\subset\gl(n,\mathbb{R})$ with
$\g^{[1]}\neq 0$ are exhausted by the subalgebras from Table
\ref{sProlC} considered as the real ones, by the subalgebras from
Table \ref{sProlR}, and subalgebras of the form
$\g_1\oplus\Real\subset\gl(n,\mathbb{R})$, where $\g_1$ is an
subalgebra from Table \ref{sProlC} considered as the real one and
with the trivial center.

{ \footnotesize
\begin{tab}\label{sProlR}  Absolutely irreducible subalgebras $\g\subset\gl(n,\Real)$ with $\g^{[1]}\neq
0$ ($\z$ denotes either $0$ or $\Real$)\\
$
\begin{array}{|c|cc|}\hline 
\g& V&\\\hline
\mathfrak{sl}(n,\Real)&\Real^n,&n\geq 3
\\\hline
\gl(n,\Real)&\Real^n,&n\geq 2
\\\hline
\mathfrak{sl}(n,\Real)\oplus\z&\odot^2\Real^n,&n\geq 3
\\\hline
\mathfrak{sl}(n,\mathbb{H})\oplus\z&S_n(\mathbb{H}),&n\geq 2
\\\hline
\mathfrak{sl}(n,\Real)\oplus\z&\Lambda^2\Real^n,&n\geq 5
\\\hline
\mathfrak{sl}(n,\mathbb{H})\oplus\z&A_n(\mathbb{H}),&n\geq 3
\\\hline
\mathfrak{sl}(n,\Real)\oplus\mathfrak{sl}(m,\Real)\oplus\Real&\Real^n\otimes\Real^m,&
n>m\geq 2\\\hline
\mathfrak{sl}(n,\Real)\oplus\mathfrak{sl}(n,\Real)\oplus\z&\Real^n\otimes\Real^n,& n\geq 3
\\\hline
\mathfrak{sl}(n,\mathbb{H})\oplus\mathfrak{sl}(m,\mathbb{H})\oplus\Real&\mathbb{H}^n\otimes\mathbb{H}^m,&
n>m\geq 1
\\\hline
\mathfrak{sl}(n,\mathbb{H})\oplus\mathfrak{sl}(n,\mathbb{H})\oplus\z&\mathbb{H}^n\otimes\mathbb{H}^n,& n\geq 2
\\\hline
\mathfrak{sl}(n,\Co)\oplus\Real&H_n(\Co),&n\geq 3\\\hline
\so(p,q)\oplus\z&\Real^{p+q},&p+q\geq 4
\\\hline
\sp(2n,\Real)\oplus\Real&\Real^{2n},&n\geq 2
\\\hline
\g\oplus\z,\,\, \g \text{ is a real form of a} &\g&\\ \text{simple complex Lie algebra}&& 
\\\hline
\end{array} $
\end{tab}
}

\section{Irreducible holonomy algebras of not symmetric odd Riemannian supermanifolds}\label{SBa}

Section \ref{secSymspaces} provides the classification of
irreducible symmetric skew-Berger algebras
$\g\subset\sp(2m,\Real)$, hence it is left to classify irreducible
non-symmetric skew-Berger algebras $\g\subset\sp(2m,\Real)$.

Irreducible skew-Berger subalgebras $\g\subset\gl(n,\Co)$ are classified  in \cite{Skew-Berger}. In Table \ref{TabsBaC} we
list irreducible skew-Berger subalgebras $\g\subset\sp(2m,\Co)$.

{ \footnotesize
\begin{tab}\label{TabsBaC} Irreducible skew-Berger subalgebras
$\g\subset\sp(2m,\Co)=\sp(V)$\\
$
\begin{array}{|c|c|c|}\hline \g& V&\text{restriction}\\\hline
\sp(2m,\Co)&\Co^{2m}&n\geq 1\\\hline
\mathfrak{sl}(2,\Co)\oplus\so(m,\Co)&\Co^2\otimes\Co^{m}&m\geq
3\\\hline \spin(12,\Co)&\Delta^+_{12}=\Co^{32}&\\\hline
\mathfrak{sl}(6,\Co)&\Lambda^3\Co^6=\Co^{20}&\\\hline
\sp(6,\Co)&V_{\pi_3}=\Co^{14}&\\\hline\hline
\so(n,\Co)\oplus\sp(2q,\Co)&\Co^n\otimes\Co^{2q}&n\geq 3,\,q\geq
2\\\hline
G^\Co_2\oplus\mathfrak{sl}(2,\Co)&\Co^7\otimes\Co^{2}&\\\hline
\so(7,\Co)\oplus\mathfrak{sl}(2,\Co)&\Co^8\otimes\Co^{2}&\\\hline
\end{array}$
\end{tab}
}

Note that for the last three subalgebras from Table \ref{TabsBaC}
it holds $\bar\R(\g)\simeq\Co$ and this space  is annihilated by
$\g$, i.e. those algebras are symmetric skew-Berger algebras.

We get now the list of irreducible not symmetric skew-Berger
subalgebras $\g\subset\sp(2m,\Real)$. Let $V$ be a real vector
space and $\g\subset\gl(V)$ an irreducible subalgebra. Consider
the complexifications $V_\Co=V\otimes_\Real\Co$ and
$\g_\Co=\g\otimes_\Real\Co\subset\gl(V_\Co)$. It is easy to see
that $\RR(\g_\Co)=\RR(\g)\otimes_\Real \Co$. Suppose that
$\g\subset\gl(V)$ is  absolutely irreducible, i.e. there exists a
complex structure $J$ on $V$ commuting with the elements of $\g$.
Then $V$ can be considered as a complex vector space. Consider the
natural representation $i:\g_\Co\to\gl(V)$ in the complex vector
space $V$. The following proposition is the analog of Proposition
3.1 from \cite{Sch}.

\begin{prop}\label{Real-Complex} Let $V$ be a real vector space and $\g\subset\gl(V)$ an irreducible subalgebra.
\begin{description}
\item[1.] If the subalgebra $\g\subset\gl(V)$ is absolutely irreducible, then $\g\subset\gl(V)$ is a skew-Berger algebra
if and only if $\g_\Co\subset\gl(V_\Co)$ is a skew-Berger algebra.

\item[2.] If the subalgebra $\g\subset\gl(V)$ is not absolutely irreducible
and if $(i(\g_\Co))^{[1]}=0$, then  $\g\subset\gl(V)$ is a skew-Berger algebra if and only if $J\g=\g$ and
$\g\subset\gl(V)$ is a complex irreducible skew-Berger algebra.
\end{description}\end{prop}

First of all, the algebras of Table \ref{TabsBaC} exhaust
the second possibility of the proposition with $(i(\g_\Co))^{[1]}=0$, and only the first 5 algebras of Table \ref{TabsBaC} are not symmetric skew-Berger algebras.

Suppose that the subalgebra $\g\subset\sp(V)$ is an absolutely
irreducible not symmetric skew-Berger algebra.  Then
$\g_\Co\subset\sp(V_\Co)$ is one of the first 5 algebras of Table
\ref{TabsBaC}. Note that each of these algebras is also a Berger
algebra \cite{Sch}. From Proposition \ref{Real-Complex}  and
Proposition 3.1 from \cite{Sch} it follows that $\g\subset\sp(V)$
is a Berger algebra, hence we may deduce all  absolutely
irreducible not symmetric skew-Berger algebra $\g\subset\sp(V)$
from \cite{Sch}.

We are left with the not absolutely irreducible subalgebras
$\g\subset\sp(V)$ such that $(i(\g_\Co))^{[1]}\neq 0$. Consider
any such $\g$. Let as in Section \ref{Prol}, $\g_1=\g\cap J\g$.
Then $\g=\g_1\oplus\g_2$ and
$\g_{\Co}=\g'_1\oplus\g''_1\oplus(\g_2\otimes \Co)$. More over,
$\g'_1$ and $\g''_1$ are isomorphic. Table \ref{sProlC} implies
that the only possible $i(\g_\Co)$ is
$\mathfrak{sl}(n,\Co)\oplus\mathfrak{sl}(n,\Co)\oplus \Co$. Then
$V=\Co^n\otimes\Co^n$, $\g_1=\mathfrak{sl}(n,\Co)$ and
$\g_2=\Real$. The Lie algebra $\g_\Co$ acting in $V_\Co$ preserves
the decomposition $V_\Co=W\oplus \bar W$, where $W$ and $\bar W$
are the eigenspaces of the extension of $J$ to $V_\Co$
corresponding to the eigenvalues $i$ and $-i$, respectively.
Moreover, $\g'_1\simeq\mathfrak{sl}(n,\Co)$ annihilates $W$,
$\g''_1\simeq\mathfrak{sl}(n,\Co)$ annihilates $\bar W$, and $\Co$
acts diagonally in  $W\oplus \bar W$. This shows that
$$\bar\R(\g_\Co)\subset\bar\R\big(\gl(n,\Co)\subset
\gl(W)\big)\oplus\bar\R\big(\gl(n,\Co)\subset \gl(\bar W)\big).$$
Note that $\dim_{\Co} W=2n$. From \cite{Skew-Berger} it follows
that the both $\mathfrak{sl}(n,\Co)$ and $\gl(n,\Co)$ do not
appear as the skew-Berger subalgebras of $\gl(W)$ for $W$ of such
dimension, hence $\bar\R(\g_\Co)=0$. This shows that $\g_1=0$,
hence $\g\cap J\g=0$, i.e. $\g\subset\g_\Co$ is a Real form. We
have to consider the real forms $\g$ of the algebras $\h$
appearing in Table \ref{sProlC} such that the restriction to $\g$
of the corresponding representation $\h\subset\gl(E)$ is
irreducible and take $V=E$ considered as the real vector space.
Now $\g_\Co=\h$ acts diagonally in $V_\Co=W\oplus \bar W$. Let
$R\in\bar\R(\g_\Co\subset V_\Co)$ and
$S\in\nabla\bar\R(\g_\Co\subset V_\Co)$.
 From the first Bianchi identity it follows that $R(X,Y)=0$ whenever both $X$ and $Y$ belong either to $W$ or to $\bar W$.
Next, for each  $X_1\in\bar W$ and $X\in W$, it holds
$$R(X_1,\cdot|_W)|_W\in((\g_\Co)_W\subset\gl(W))^{[1]},\quad
R(X,\cdot|_{\bar W})|_{\bar W}\in((\g_\Co)_{\bar W}\subset\gl(\bar
W))^{[1]}.$$ From this and the second Bianchi identity it follows
that
$$S_{\cdot|_W}(X_1,\cdot|_W)|_W\in((\g_\Co)_W\subset\gl(W))^{[2]},\quad
S_{\cdot|_{\bar W}}(X,\cdot|_{\bar W})|_{\bar W}\in((\g_\Co)_{\bar
W}\subset\gl(\bar W))^{[2]}.$$ Hence, if $\g\subset\sp(V)$ is a
skew-Berger algebra and $((\g_\Co)_W\subset\gl(W))^{[2]}=0$, then
$\nabla\bar\R(\g)=0$, i.e. $\g\subset\sp(V)$ is a symmetric
skew-Berger algebra. Thus we get only the following  4
possibilities for $\h\subset\gl(E)$:
$$\gl(n,\Co),\ \mathfrak{sl}(n,\Co),\ \so(n,\Co),\ \so(n,\Co)\oplus\Co \subset\gl(n,\Co).$$
The corresponding $\g\subset\sp(2m,\Real)$ are the following:
$$\u(n),\ \su(n)\subset\sp(2n,\Real),\quad \so(n,\H),\ \so(n,\H)\oplus\Real i\subset\sp(4n,\Real).$$
Let us find the spaces
$\bar\R(\g)\otimes\Co=\bar\R(\g_\Co\subset\gl(V_\Co))$. As we have
seen above, any $R\in \bar \R(\g_\Co\subset\gl(V_\Co))$ is
uniquely defined by the values $R(X,X_1)$, where $X\in W$ and
$X_1\in\bar W$. Let $e_1,...,e_m$ be a basis in $W$ and
$e^1,...,e^m$ the dual basis in $\bar W$. We may write
$R(e_i,e^j)=A^j_i$ for some $A^j_i\in g_{\Co}|_{W}$. Define the
numbers $A_{ik}^{jl}$ such that $A^j_i e_k=\sum_l A_{ik}^{jl}e_l$.
Then $A^j_i e^l=\sum_k A_{ik}^{jl}e^k$. Let $\g=\u(n)$, then
$A^j_i\in\gl(n,\Co)$ and $m=n$. From above it follows that
$R\in\bar R(\g_\Co)$ if and only if $A_{ik}^{jl}=-A_{ki}^{jl}$ and
$A_{ik}^{jl}=-A_{ik}^{lj}$. Hence, $\bar\R(\u(n))\otimes\Co$ is
isomorphic $(W\wedge W)\otimes(W^*\wedge W^*)$. For  $\g=\su(n)$
we get the additional condition $\sum_k A_{ik}^{jk}=0$. This shows
that $\u(n),\ \su(n)\subset\sp(2n,\Real)$ are skew-Berger
subalgebras. Since
$$(\so(2n,\Co)\subset\gl(2n,\Co))^{[1]}=(\so(2n,\Co)\oplus\Co\subset\gl(2n,\Co))^{[1]},$$
$\so(n,\H)\oplus\Real i \subset\sp(4n,\Real)$ is not a skew-Berger
subalgebra. Finally consider $\g=\so(n,\H) \subset\sp(4n,\Real)$.
Any $R\in\bar\R(\g_\Co\subset V_\Co)$ is defined as above by the
numbers $A_i^j$ with the additional condition
$A_{ik}^{jl}=-A_{il}^{jk}$. This shows that $\bar\R(\g_\Co\subset
V_\Co)\simeq\wedge^4\Co^{2n}$, i.e. $\so(n,\H)
\subset\sp(4n,\Real)$ is a skew-Berger subalgebra. We obtain the
following classification theorem

\begin{theorem}  Possible irreducible holonomy algebras
$\g\subset\sp(2m,\Real)$ of not symmetric odd Riemannian
supermanifolds are listed in Table \ref{TabsBaR}.\end{theorem}

{\footnotesize
\begin{tab}\label{TabsBaR} Possible irreducible holonomy algebras
$\g\subset\sp(2m,\Real)=\sp(V)$ of not symmetric odd Riemannian
supermanifolds.\\
$
\begin{array}{|c|c|c|}\hline \g& V&\text{restriction}\\\hline
\sp(2m,\Real)&\Real^{2m}&m\geq 1\\\hline
 \u(p,q)&\mathbb{C}^{p,q}&p+q\geq
2\\\hline \su(p,q)&\mathbb{C}^{p,q}&p+q\geq 2\\\hline
\so(n,\mathbb{H})&\mathbb{H}^{n}&n\geq
2\\\hline\sp(1)\oplus\so(n,\mathbb{H})&\mathbb{H}^{n}&n\geq
2\\\hline
\mathfrak{sl}(2,\Real)\oplus\so(p,q)&\Real^2\otimes\Real^{p,q}&p+q\geq
3\\\hline  \spin(2,10)&\Delta^+_{2,10}=\Real^{32}&\\\hline
\spin(6,6)&\Delta^+_{6,6}=\Real^{32}&\\ \hline
\so(6,\mathbb{H})&\Delta^\mathbb{H}_{6}=\mathbb{H}^{8}&\\\hline
\mathfrak{sl}(6,\Real)&\Lambda^3\Real^6=\Real^{20}&\\\hline
\su(1,5)&\{\omega\in\Lambda^3\Co^6|*w=w\}&\\\hline
\su(3,3)&\{\omega\in\Lambda^3\Co^6|*w=w\}&\\\hline
\sp(6,\Real)&\Real^{14}\subset\Lambda^3\Real^6&\\\hline\sp(2m,\Co)&\Co^{2m}&m\geq
1\\\hline
\mathfrak{sl}(2,\Co)\oplus\so(m,\Co)&\Co^2\otimes\Co^{m}&m\geq
3\\\hline \spin(12,\Co)&\Delta^+_{12}=\Co^{32}&\\\hline
\mathfrak{sl}(6,\Co)&\Lambda^3\Co^6=\Co^{20}&\\\hline
\sp(6,\Co)&V_{\pi_3}=\Co^{14}&\\\hline
\end{array}$
\end{tab}
}

\bibliographystyle{unsrt}

\end{document}